# UNIFORM IN BANDWIDTH CONSISTENCY OF KERNEL-TYPE FUNCTION ESTIMATORS

By Uwe Einmahl[1] and David M. Mason[2]

*Vrije Universiteit Brussel and University of Delaware*

We introduce a general method to prove uniform in bandwidth consistency of kernel-type function estimators. Examples include the kernel density estimator, the Nadaraya–Watson regression estimator and the conditional empirical process. Our results may be useful to establish uniform consistency of data-driven bandwidth kernel-type function estimators.

**1. Introduction and statements of main results.** Let $X, X_1, X_2, \ldots$ be i.i.d. $\mathbb{R}^d$, $d \geq 1$, valued random variables and assume that the common distribution function of these variables has a Lebesgue density function, which we shall denote by $f$. A kernel $K$ will be any measurable function which satisfies the conditions

(K.i) $$\int_{\mathbb{R}^d} K(s) \, ds = 1,$$

(K.ii) $$\|K\|_{\infty} := \sup_{x \in \mathbb{R}^d} |K(x)| = \kappa < \infty.$$

The kernel density estimator of $f$ based upon the sample $X_1, \ldots, X_n$ and bandwidth $0 < h < 1$ is

$$\widehat{f}_{n,h}(x) = (nh)^{-1} \sum_{i=1}^{n} K((x - X_i)/h^{1/d}), \qquad x \in \mathbb{R}^d.$$

Choosing a suitable bandwidth sequence $h_n \to 0$ and assuming that the density $f$ is continuous, one obtains a strongly consistent estimator $\widehat{f}_n := \widehat{f}_{n,h_n}$

Received March 2003; revised April 2004.

[1]Supported in part by an FWO grant.

[2]Supported in part by NSA Grant MDA904-02-1-0034 and NSF Grant DMS-02-03865.

of $f$, that is, one has with probability 1, $\widehat{f}_n(x) \to f(x)$, $x \in \mathbb{R}^d$. There are also results concerning uniform convergence and convergence rates. For proving such results one usually writes the difference $\widehat{f}_n(x) - f(x)$ as the sum of a probabilistic term $\widehat{f}_n(x) - \mathbb{E}\widehat{f}_n(x)$ and a deterministic term $\mathbb{E}\widehat{f}_n(x) - f(x)$, the so-called bias. The order of the bias depends on smoothness properties of $f$ only, whereas the first (random) term can be studied via empirical process techniques, as has been pointed out by Stute [29, 30, 31] and Pollard [26], among other authors.

A recent result by Giné and Guillou [14] (see also [5]) shows that if $K$ is a "regular" kernel, the density function $f$ is bounded and $h_n$ satisfies the regularity conditions $h_n \searrow 0$, $h_n/h_{2n}$ is bounded,

$$\log(1/h_n)/\log\log n \to \infty \quad \text{and} \quad nh_n/\log n \to \infty,$$

one has with probability 1,

(1.1)               $\|\widehat{f}_n - \mathbb{E}\widehat{f}_n\|_\infty = O(\sqrt{\log(1/h_n)/nh_n})$.

Moreover, this rate cannot be improved. Interestingly, one does not need continuity of $f$ for this result. (Of course, continuity of $f$ is crucial for controlling the bias.)

Some related results on uniform convergence over compact subsets have been obtained by Einmahl and Mason (EM) [11] for a much larger class of estimators including kernel estimators for regression functions among others. In this general setting, however, it is often not possible to obtain the convergence uniformly over $\mathbb{R}^d$. Density estimators are in that sense somewhat exceptional.

The main purpose of this paper is to introduce a method to establish consistency of kernel-type estimators when the bandwidth $h$ is allowed to range in a small interval which may decrease in length with the sample size. Our results will be immediately applicable to proving uniform consistency of kernel-type estimators when the bandwidth $h$ is a function of the location $x$ or the data $X_1, \ldots, X_n$. The resulting "variable bandwidth kernel estimators" are from a statistical point of view clearly preferable to those bandwidths which are only a function of the sample size $n$, ignoring the data and the location. We discuss this in more detail in Remark 7 below, after we have stated some of our main results. Furthermore, we address the issue of bias in Remark 6.

In order to formulate our results let us first specify what we mean by a "regular" kernel $K$. Consider the class of functions

$$\mathcal{K} = \{K((x-\cdot)/h^{1/d}) : h > 0, x \in \mathbb{R}^d\}.$$

For $\varepsilon > 0$, let $N(\varepsilon, \mathcal{K}) = \sup_Q N(\kappa\varepsilon, \mathcal{K}, d_Q)$, where the supremum is taken over all probability measures $Q$ on $(\mathbb{R}^d, \mathcal{B})$, $d_Q$ is the $L_2(Q)$-metric and, as usual, $N(\varepsilon, \mathcal{K}, d_Q)$ is the minimal number of balls $\{g : d_Q(g, g') < \varepsilon\}$ of $d_Q$-radius $\varepsilon$ needed to cover $\mathcal{K}$. Assume that $\mathcal{K}$ satisfies the following uniform entropy condition:

(K.iii) for some $C > 0$ and $\nu > 0$, $N(\varepsilon, \mathcal{K}) \leq C\varepsilon^{-\nu}$, $0 < \varepsilon < 1$.

Pollard [26], Nolan and Pollard [25] and van der Vaart and Wellner [35] provide a number of sufficient conditions for (K.iii) to hold. For instance, it is satisfied for general $d \geq 1$ whenever $K(x) = \phi(p(x))$, with $p(x)$ being a polynomial in $d$ variables and $\phi$ being a real-valued function of bounded variation.

Finally, to avoid using outer probability measures in all of our statements, we impose the following measurability assumption.

(K.iv) $\mathcal{K}$ is a pointwise measurable class, that is, there exists a countable subclass $\mathcal{K}_0$ of $\mathcal{K}$ such that we can find for any function $g \in \mathcal{K}$ a sequence of functions $\{g_m\}$ in $\mathcal{K}_0$ for which

$$g_m(z) \to g(z), \qquad z \in \mathbb{R}^d.$$

This condition is discussed in [35]. It is satisfied whenever $K$ is right continuous.

Our first result concerning density estimators is the following.

THEOREM 1. *Assuming* (K.i)–(K.iv) *and* $f$ *is bounded, we have for any* $c > 0$, *with probability* 1,

$$(1.2) \qquad \limsup_{n \to \infty} \sup_{c \log n/n \leq h \leq 1} \frac{\sqrt{nh} \|\widehat{f}_{n,h} - \mathbb{E}\widehat{f}_{n,h}\|_\infty}{\sqrt{\log(1/h) \vee \log\log n}} =: K(c) < \infty.$$

REMARK 1. Though this was not our main goal, we point out that if one chooses a deterministic sequence $h_n$ satisfying $nh_n/\log n \to \infty$ and $\log(1/h_n)/\log\log n \to \infty$, one re-obtains (1.1), which is Theorem 1 of Giné and Guillou [14] with slightly less regularity. (We do not need to assume, as they do, that $h_n \searrow 0$ or that $h_n/h_{2n}$ is bounded.)

REMARK 2. With applications to variable bandwidth estimators in mind, we further note that Theorem 1 implies for any sequences $0 < a_n < b_n \leq 1$, satisfying $b_n \to 0$ and $na_n/\log n \to \infty$, with probability 1,

$$(1.3) \qquad \sup_{a_n \leq h \leq b_n} \|\widehat{f}_{n,h} - \mathbb{E}\widehat{f}_{n,h}\|_\infty = O\left(\sqrt{\frac{\log(1/a_n) \vee \log\log n}{na_n}}\right),$$

which in turn implies

$$(1.4) \qquad \lim_{n \to \infty} \sup_{a_n \leq h \leq b_n} \|\widehat{f}_{n,h} - \mathbb{E}\widehat{f}_{n,h}\|_\infty = 0 \qquad \text{a.s.}$$

REMARK 3. It is routine to modify the proof of Theorem 1 to show that it remains true when (K.iii) is replaced by the bracketing condition:

(K′.iii) for some $C_0 > 0$ and $\nu_0 > 0$, $N_{[\cdot]}(\varepsilon, \mathcal{F}, L_2(P)) \leq C_0 \varepsilon^{-\nu_0}$, $0 < \varepsilon < 1$.

Refer to page 270 of [34] for the definition of $N_{[\cdot]}(\varepsilon, \mathcal{F}, L_2(P))$. Essentially all that one has to do is to substitute the use of Corollary 4 by Lemma 19.34 of van der Vaart [34].

For a related result refer to Theorem 1 of Nolan and Marron [24], where almost sure convergence to zero has been established in a similar setting. On the other hand, our result provides explicit convergence rates for kernel density estimators.

Let us now look at the bias term. As soon as we know that

$$\sup_{a_n \leq h \leq b_n} \|\mathbb{E}\widehat{f}_{n,h} - f\|_\infty \to 0, \tag{1.5}$$

we have under the conditions of Theorem 1,

$$\sup_{a_n \leq h \leq b_n} \|\widehat{f}_{n,h} - f\|_\infty \to 0.$$

If $f$ is uniformly continuous on $\mathbb{R}^d$, here is a sufficient condition for (1.5) which is easy to verify: Define

$$\Psi_K(x) = \sup_{|y| \geq |x|} |K(y)|, \qquad x \in \mathbb{R}^d,$$

and introduce the assumption

$$\int_{\mathbb{R}^d} \Psi_K(x)\, dx < \infty. \tag{K.v}$$

Note that this assumption trivially holds for a compactly supported kernel function.

COROLLARY 1. Assuming (K.i)–(K.v) for any sequences $0 < a_n < b_n < 1$, satisfying $b_n \to 0$ and $na_n/\log n \to \infty$, and any uniformly continuous density $f$, we have

$$\lim_{n \to \infty} \sup_{a_n \leq h \leq b_n} \|\widehat{f}_{n,h} - f\|_\infty = 0 \qquad a.s. \tag{1.6}$$

REMARK 4. If $a_n = c \log n/n$ for some $c > 0$, then (1.6) does not hold, that is, the limit in (1.6) is positive. Refer to [4] and [6] for details.

Our method is not restricted to the case of kernel density estimators. To give the reader an indication of what other kinds of kernel-type estimators can be treated using our techniques, consider i.i.d. $(d+1)$-dimensional random vectors $(Y, X)$, $(Y_1, X_1)$, $(Y_2, X_2), \ldots$, where the $Y$-variables are one-dimensional. We shall assume that $X$ has a marginal Lebesgue density function $f$ and that the regression function

$$m(x) = \mathbb{E}[Y|X=x], \qquad x \in \mathbb{R}^d,$$

exists. Let $\widehat{m}_{n,h}(x)$ be the usual Nadaraya–Watson estimator of $m(x)$ with bandwidth $0 < h < 1$, that is,

$$\widehat{m}_{n,h}(x) = \frac{\sum_{i=1}^n Y_i K((x - X_i)/h^{1/d})}{\sum_{i=1}^n K((x - X_i)/h^{1/d})}.$$

A huge literature has been developed on the consistency of the Nadaraya–Watson estimator. Consult [16] and [11] for references to some of the more important work.

Assuming that $m$ is $p+1$ times differentiable at a fixed $x_0$, one can use the local polynomial regression techniques of Fan and Gijbels [12] to obtain a better estimate at $x_0$ than that given by the Nadaraya–Watson estimator. We will not treat the uniform consistency of such estimators in the present paper. It should, however, be feasible to apply similar empirical process methods in this setting as well.

With the above setup we have the following uniform in bandwidth result. Set

$$\bar{r}(x, h) = \mathbb{E}[Y K((x - X)/h^{1/d})]/h \quad \text{and} \quad \bar{f}(x, h) = \mathbb{E}[K((x - X)/h^{1/d})]/h.$$

For any subset $I$ of $\mathbb{R}^d$, let $I^\varepsilon$ denote its closed $\varepsilon$-neighborhood with respect to the maximum-norm $|\cdot|_+$ on $\mathbb{R}^d$, that is, $|x|_+ = \max_{1 \le i \le d} |x_i|$, $x \in \mathbb{R}^d$. Set further for any function $\psi : \mathbb{R}^d \to \mathbb{R}$, $\|\psi\|_I = \sup_{x \in I} |\psi(x)|$.

THEOREM 2. *Let $I$ be a compact subset of $\mathbb{R}^d$ and let $K$ satisfy* (K.i)–(K.iv) *with support contained in* $[-1/2, 1/2]^d$. *Suppose further that there exists an $\varepsilon > 0$ so that*

(1.7)      *$f$ is continuous and strictly positive on $J := I^\varepsilon$.*

*If there exists an $M > 0$ such that*

(1.8)      $|Y| \mathbb{1}\{X \in J\} \le M \qquad a.s.,$

*we have for large enough $c > 0$ and any $b_n \searrow 0$,*

(1.9)      $\displaystyle \limsup_{n \to \infty} \sup_{c \log n/n \le h \le b_n} \frac{\sqrt{nh} \|\widehat{m}_{n,h} - \bar{r}(\cdot, h)/\bar{f}(\cdot, h)\|_I}{\sqrt{\log(1/h) \vee \log \log n}}$

         $=: K(I, c) < \infty \qquad a.s.$

*Moreover, if instead of* (1.8) *we assume that for some* $p > 2$

$$(1.10) \qquad \sup_{z \in J} \mathbb{E}(|Y|^p | X = z) =: \alpha < \infty,$$

*we have for any* $c > 0$ *and* $b_n \searrow 0$ *with* $\gamma = \gamma(p) = 1 - 2/p,$

$$
(1.11)
\begin{aligned}
&\limsup_{n \to \infty} \sup_{c(\log n/n)^\gamma \le h \le b_n} \frac{\sqrt{nh} \|\widehat{m}_{n,h} - \bar{r}(\cdot, h)/\bar{f}(\cdot, h)\|_I}{\sqrt{\log(1/h) \vee \log \log n}} \\
&\qquad =: K'(I, c) < \infty \qquad a.s.
\end{aligned}
$$

COROLLARY 2. *Let* $I$ *be a compact subset of* $\mathbb{R}^d$ *and let* $K$ *satisfy* (K.i)–(K.iv) *with support contained in* $[-1/2, 1/2]^d$. *Assume that the distribution function of* $(Y, X)$ *has a Lebesgue density* $(y, x) \to p(y, x)$, *so that the marginal density of* $X$ *is given by*

$$f(x) = \int_{-\infty}^{\infty} p(y, x) \, dy, \qquad x \in \mathbb{R}^d.$$

*Suppose further that there exists an* $\varepsilon > 0$ *so that* (1.7) *holds and that*

$$(1.12) \quad \text{for all } z \in J, \qquad \lim_{z' \to z} p(y, z') = p(y, z) \qquad \text{for almost every } y \in \mathbb{R}.$$

*If* (1.8) *holds, then for* $0 < a_n < b_n < 1$, *satisfying* $b_n \to 0$ *and* $na_n/\log n \to \infty$,

$$(1.13) \qquad \lim_{n \to \infty} \sup_{a_n \le h \le b_n} \|\widehat{m}_{n,h} - m(\cdot)\|_I = 0 \qquad a.s.$$

*If* (1.10) *holds, then with* $\gamma = 1 - 2/p$ *for* $0 < c(\log n/n)^\gamma < b_n < 1$ *satisfying* $b_n \to 0$,

$$(1.14) \qquad \lim_{n \to \infty} \sup_{c(\log n/n)^\gamma \le h \le b_n} \|\widehat{m}_{n,h} - m(\cdot)\|_I = 0 \qquad a.s.$$

REMARK 5. Let us also mention that if, in the bounded case, we choose a deterministic bandwidth sequence $h_n$ satisfying the standard assumption $nh_n/\log n \to \infty$ and $\log(1/h_n)/\log \log n \to \infty$, we get that with probability 1,

$$\limsup_{n \to \infty} \frac{\sqrt{nh_n} \|\widehat{m}_{n,h_n} - \bar{r}(\cdot, h_n)/\bar{f}(\cdot, h_n)\|_I}{\sqrt{2 \log(1/h_n)}} \le C < \infty.$$

This is a sharp result. In our previous paper [11] we have shown under additional assumptions ($h_n \searrow 0$ and $nh_n \nearrow$, $d = 1$, $I = [a, b]$ and $K$ satisfies a continuity condition and is of bounded variation on $\mathbb{R}$) that the lim sup is positive and actually a limit. [Note, however, that the limiting constant has not been correctly stated in formula (1.16) of that paper. With the notation of the present paper the limiting constant is $\sup_{x \in I} (\sigma(x) \|K\|_2)/\sqrt{f(x)}$, where $\sigma^2(x) = \text{Var}(Y|X = x)$.] Moreover, if (1.8) holds, then a result of Collomb [3] implies that the condition $na_n/\log n \to \infty$ is necessary for uniform consistency.

REMARK 6. Under additional smoothness assumptions on $f$ one can also derive explicit convergence rates in (1.6) and (1.14). For instance, if one knows that $f$ is uniformly Lipschitz continuous, one easily sees that the bias (1.5) is of order $O(b_n^{1/d})$, which permits one to derive a convergence rate in (1.6) one which depends on $a_n$, via the rate from Theorem 1, and on $b_n$, via the rate in (1.5). For more information on the interplay between smoothness and the size of the bias term consult [1, 8, 10]. Similarly under extra smoothness conditions the bias term in the Nadaraya–Watson estimator is well behaved and one also can specify convergence rates. For appropriate smoothness conditions refer to [1] and, especially, to Section 2.3 of [7].

REMARK 7. Suppose now that $\widehat{h}_n = \widehat{h}_n(x)$ is a local data-driven bandwidth sequence satisfying

$$(1.15) \qquad \mathbb{P}\{a_n \le \widehat{h}_n(x) \le b_n : x \in I\} \to 1,$$

or a constant data-driven bandwidth sequence $\widehat{h}_n$ satisfying with probability 1, for all large enough $n \ge 1$,

$$(1.16) \qquad a_n \le \widehat{h}_n \le b_n.$$

For instance, if $d = 1$, one often has for appropriate $0 < a < b < \infty$, $a_n = an^{-1/5}$ and $b_n = bn^{-1/5}$. [10] is a good place to read about the various optimality criteria that lead to the $n^{-1/5}$. In this case and more generally under the assumptions of Corollary 1,

$$\|\widehat{f}_{n,\widehat{h}_n} - f\|_I \to 0,$$

and under those of Corollary 2,

$$\|\widehat{m}_{n,\widehat{h}_n} - m(\cdot)\|_I \to 0,$$

where the convergence is either in probability or with probability 1 depending on whether (1.15) or (1.16) holds.

Deheuvels and Mason [7] consider local plug-in type estimators $\hat{h}_n(x)$ which satisfy (1.15) with $a_n = c_1 h_n$ and $b_n = c_2 h_n$, where $c_1 < c_2$, or

$$(1.17) \qquad \mathbb{P}\left(\sup_{x \in I} |\hat{h}_n(x)/h_n - C(x)| > \varepsilon\right) \to 0$$

for any $\varepsilon > 0$, where $C$ is an appropriate continuous function on $I$. Refer especially to their Example 2.1, where they show subject to smoothness assumptions that the optimal $\hat{h}_n(x)$ in terms of asymptotic mean square error for estimating $f$ or $m$ satisfies (1.17) with $h_n = n^{-1/5}$.

The literature on data-driven bandwidth selection is extensive. We cite, for instance, [[2, 17, 21, 22, 23, 27]]. For further references and methods consult [18], Chapter 7 of [10], [7] and [9].

All data-driven bandwidth selection procedures require some smoothness assumptions in order to get rates. Our results show that even if such assumptions do not hold, one may still have consistency as long as (1.15) is satisfied for appropriate $a_n$ and $b_n$ not necessarily of the form $a_n = c_1 h_n$ and $b_n = c_2 h_n$.

Our next example is a kernel estimator of the conditional distribution function

$$F(t|z) := \mathbb{P}(Y \le t | X = z),$$

defined for a kernel $K$ and bandwidth $0 < h < 1$ to be

$$(1.18) \qquad \widehat{F}_{n,h}(t|z) := \frac{\sum_{i=1}^n \mathbb{1}(Y_i \le t) K((z - X_i)/h^{1/d})}{\sum_{i=1}^n K((z - X_i)/h^{1/d})}.$$

Stute [32] calls this the *conditional empirical distribution function* and was the first to establish uniform consistency results for it.

THEOREM 3. *Let $I$ be a compact subset of $\mathbb{R}^d$ and let $K$ satisfy* (K.i)–(K.iv) *with support contained in $[-1/2, 1/2]^d$. Suppose further that there exists an $\varepsilon > 0$ so that* (1.7) *holds. Then, with probability* 1, *we have for large enough $c > 0$ and any $b_n \searrow 0$,*

$$(1.19) \qquad \limsup_{n \to \infty} \sup_{c \log n / n \le h \le b_n} \frac{\sup_{z \in I} \sqrt{nh}\, \|\widehat{F}_{n,h}(\cdot|z) - F_{n,h}(\cdot|z)\|_\infty}{\sqrt{\log(1/h) \vee \log\log n}}$$

$$=: K''(I, c) < \infty,$$

*where $F_{n,h}(t|z) = \mathbb{E}[K((z - X)/h^{1/d}) \mathbb{1}\{Y \le t\}]/(h \mathbb{E}\widehat{f}_{n,h}(z))$, $t \in \mathbb{R}$.*

COROLLARY 3. *Let $I$ be a compact subset of $\mathbb{R}^d$ and let $K$ satisfy* (K.i)–(K.iv) *with support contained in $[-1/2, 1/2]^d$. Suppose further that there exists an $\varepsilon > 0$ so that* (1.7) *holds and* (1.12) *is satisfied. Then for $0 < a_n < b_n < 1$, satisfying $b_n \to 0$ and $na_n/\log n \to \infty$,*

$$(1.20) \qquad \lim_{n \to \infty} \sup_{a_n \le h \le b_n} \sup_{z \in I} \|\widehat{F}_{n,h}(\cdot|z) - F(\cdot|z)\|_\infty = 0.$$

REMARK 8. Sometimes one wants to use vector bandwidths (see, in particular, Chapter 12 of Devroye and Lugosi [9]). With obvious changes of notation, our results and their proofs remain true when $h_n$ is replaced by a vector bandwidth $\mathbf{h}_n = (h_n^{(1)}, \ldots, h_n^{(d)})$, where $\min_{1 \leq i \leq d} h_n^{(i)} > 0$. In this situation we set $h_n = \prod_{i=1}^d h_n^{(i)}$, and for any vector $\mathbf{v} = (v_1, \ldots, v_d)$ we replace $\mathbf{v}/h_n^{1/d}$ by $(v_1/h_n^{(1)}, \ldots, v_d/h_n^{(d)})$. For ease of presentation we chose to use real-valued bandwidths throughout.

Theorem 1 is proved in Section 2. Theorems 2 and 3 will follow from a more general result stated and proved in Section 3. Our proofs are based on an extension of the methods developed in [11]. We use the same idea which was developed in [11], namely, combining an exponential inequality of Talagrand [33] with a suitable moment inequality.

## 2. Proofs of Theorem 1 and Corollary 1.
We shall look at a slightly more general setup than in the Introduction. Let $(\mathcal{X}, \mathcal{A})$ be a measurable space. Throughout this section we assume that on our basic probability space $(\Omega, \mathcal{F}, \mathbb{P})$ we have independent $(\mathcal{F}, \mathcal{A})$-measurable variables $X_i : \Omega \to \mathcal{X}$, $1 \leq i \leq n$, with common distribution $\mu$.

Let $\mathcal{G}$ be a pointwise measurable class of functions from $\mathcal{X}$ to $\mathbb{R}$ (see the Introduction and Example 2.3.4 in [35]). Further let $\varepsilon_1, \ldots, \varepsilon_n$ be a sequence of independent Rademacher random variables, independent of $X_1, \ldots, X_n$. Let $G$ be a finite-valued measurable function satisfying for all $x \in \mathcal{X}$,

$$(2.1) \qquad G(x) \geq \sup_{g \in \mathcal{G}} |g(x)|,$$

and define

$$(2.2) \qquad N(\varepsilon, \mathcal{G}) = \sup_Q N(\varepsilon \sqrt{Q(G^2)}, \mathcal{G}, d_Q),$$

where the supremum is taken over all probability measures $Q$ on $(\mathcal{X}, \mathcal{A})$ for which $0 < Q(G^2) < \infty$ and $d_Q$ is the $L_2(Q)$-metric.

We need the following version of Proposition A.1 of EM [11].

PROPOSITION 1. *Let $\mathcal{G}$ be a pointwise measurable class of bounded functions such that for some constants $C, \nu \geq 1$ and $0 < \sigma \leq \beta$ and $G$ as above, the following conditions hold:*

(i) $\mathbb{E}[G(X)^2] \leq \beta^2$;

(ii) $N(\varepsilon, \mathcal{G}) \leq C\varepsilon^{-\nu}, 0 < \varepsilon < 1$;

(iii) $\sigma_0^2 := \sup_{g \in \mathcal{G}} \mathbb{E}[g(X)^2] \leq \sigma^2$;

(iv) $\sup_{g \in \mathcal{G}} \|g\|_\infty \leq \frac{1}{4\sqrt{\nu}} \sqrt{n\sigma^2/\log(C_1\beta/\sigma)}$, *where $C_1 = C^{1/\nu} \vee e$.*

*Then we have for some absolute constant $A$,*

$$(2.3) \qquad \mathbb{E}\left\|\sum_{i=1}^{n}\varepsilon_i g(X_i)\right\|_{\mathcal{G}} \leq A\sqrt{\nu n\sigma^2 \log(C_1\beta/\sigma)}.$$

PROOF. Our proof is a modification of that of Proposition A.1 of EM [11]. We denote vectors $(x_1,\ldots,x_n) \in \mathcal{X}^n$ by $\mathbf{x}$ and we define the subsets $F_n$ and $G_n$ of $\mathcal{X}^n$ as in this paper, that is,

$$G_n := \left\{ \mathbf{x}: n^{-1}\sum_{j=1}^{n}G^2(x_j) \leq 256\beta^2 \right\},$$

$$F_n := \left\{ \mathbf{x}: n^{-1}\sup_{g \in \mathcal{G}}\sum_{j=1}^{n}g^2(x_j) \leq 64\sigma^2 \right\}.$$

We can infer from (A.8)–(A.10) in [11] that on $F_n \cap G_n$,

$$\mathbb{E}\left\|\sum_{i=1}^{n}\varepsilon_i g(x_i)\right\|_{\mathcal{G}} \leq K'\sigma\sqrt{n\nu \log(C_1\beta/\sigma)},$$

where $K'$ is an absolute constant. Therefore, we have for

$$t \geq 96K'\sigma\sqrt{n\nu \log(C_1\beta/\sigma)},$$

$$(2.4) \qquad \mathbb{P}\left\{\left\|\sum_{i=1}^{n}\varepsilon_i g(x_i)\right\|_{\mathcal{G}} > t\right\} \leq 1/96 \qquad \forall \mathbf{x} \in F_n \cap G_n,$$

and, consequently, in this range of $t$ that

$$\mathbb{P}\left\{\left\|\sum_{i=1}^{n}\varepsilon_i g(X_i)\right\|_{\mathcal{G}} > t\right\} \leq 1/96 + \mu^n(F_n^c) + \mu^n(G_n^c).$$

By Markov's inequality we trivially have $\mu^n(G_n^c) \leq 1/256$. Using Lemma 5.1 of [15] exactly as in [11] and recalling that $C_1\beta/\sigma \geq e$ and $\nu \geq 1$, we see that

$$\mu^n(F_n^c) \leq 4\mu^n(G_n^c) + 12 \cdot 16^{\nu}(C_1\beta/\sigma)^{-15\nu} \leq 7/256,$$

which finally implies that

$$\mathbb{P}\left\{\left\|\sum_{i=1}^{n}\varepsilon_i g(X_i)\right\|_{\mathcal{G}} > t\right\} \leq 1/24,$$

whenever (2.4) holds. A straightforward application of the Hoffmann–Jørgensen inequality as stated in Proposition 6.8 of Ledoux and Talagrand [20] finally yields the desired moment inequality. $\square$

From the above moment inequality we can infer the following:

COROLLARY 4. *Let $\mathcal{G}$ be as in Proposition* 1 *satisfying* (i)–(iii), *and instead of* (iv) *assume that*

(v) $\sup_{g \in \mathcal{G}} \|g\|_\infty \leq U$, *where* $\sigma_0 \leq U \leq C_2 \sqrt{n}\beta$, *and* $C_2 = \frac{1}{4\sqrt{\nu \log C_1}}$.

*Then we have*

$$(2.5) \quad \mathbb{E}\left\|\sum_{i=1}^n \varepsilon_i g(X_i)\right\|_{\mathcal{G}} \leq A\{\sqrt{\nu n \sigma_0^2 \log(C_1\beta/\sigma_0)} + 2\nu U \log(C_3 n(\beta/U)^2)\},$$

*where* $C_3 = C_1^2/16\nu$.

PROOF. Whenever

$$U \leq \frac{1}{4\sqrt{\nu}}\sqrt{n\sigma_0^2/\log(C_1\beta/\sigma_0)},$$

inequality (2.5) follows immediately from our proposition by choosing $\sigma = \sigma_0$.

Assume now that

$$\frac{1}{4\sqrt{\nu}}\sqrt{n\sigma_0^2/\log(C_1\beta/\sigma_0)} < U \leq C_2\sqrt{n}\beta.$$

Then using the monotonicity of the function $t \to \sqrt{nt^2/\log(C_1\beta/t)}$ we can find a unique $\sigma \in \,]\sigma_0, \beta]$ satisfying

$$U = \frac{1}{4\sqrt{\nu}}\sqrt{n\sigma^2/\log(C_1\beta/\sigma)}.$$

Applying our proposition with this choice of $\sigma$, it follows that

$$\mathbb{E}\left\|\sum_{i=1}^n \varepsilon_i g(X_i)\right\|_{\mathcal{G}} \leq A\sqrt{\nu n \sigma^2 \log(C_1\beta/\sigma)} \leq 4A\nu U \log(C_1\beta/\sigma).$$

Next rewriting the equation which defines $\sigma$ and recalling that $C_1\beta/\sigma \geq e$, we readily obtain that

$$1/\sigma \leq \sigma^{-1}\sqrt{\log(C_1\beta/\sigma)} = \sqrt{n}/(4\sqrt{\nu}U),$$

and thus

$$C_1(\beta/\sigma) \leq C_1/(4\sqrt{\nu})(\sqrt{n}\beta/U) =: \sqrt{C_3}(\sqrt{n}\beta/U).$$

It follows that

$$\mathbb{E}\left\|\sum_{i=1}^n \varepsilon_i g(X_i)\right\|_{\mathcal{G}} \leq 2A\nu U \log(C_3 n(\beta/U)^2),$$

which proves the corollary.   □

A bound similar to that given in Corollary 4 has been given by Giné and Guillou [13] using a different method.

As already indicated in the Introduction, our proof is based on an inequality of Talagrand [33] (see also [19]) which we state here for easy reference later on.

Let $\alpha_n$ be the empirical process based on the sample $X_1, \ldots, X_n$, that is, if $g: \mathcal{X} \to \mathbb{R}$, we have

$$\alpha_n(g) = \sum_{i=1}^{n} (g(X_i) - \mathbb{E}g(X))/\sqrt{n},$$

and set for any class $\mathcal{G}$ of such functions

$$\|\sqrt{n}\alpha_n\|_{\mathcal{G}} = \sup_{g \in \mathcal{G}} |\sqrt{n}\alpha_n(g)|.$$

INEQUALITY. *Let $\mathcal{G}$ be a pointwise measurable class of functions satisfying for some $0 < M < \infty$,*

$$\|g\|_{\infty} \leq M, \qquad g \in \mathcal{G}.$$

*Then we have for all $t > 0$,*

$$\mathbb{P}\left\{ \max_{1 \leq m \leq n} \|\sqrt{m}\alpha_m\|_{\mathcal{G}} \geq A_1 \left( \mathbb{E}\left\| \sum_{i=1}^{n} \varepsilon_i g(X_i) \right\|_{\mathcal{G}} + t \right) \right\}$$

$$\leq 2\left\{ \exp\left( -\frac{A_2 t^2}{n\sigma_{\mathcal{G}}^2} \right) + \exp\left( -\frac{A_2 t}{M} \right) \right\},$$

*where $\sigma_{\mathcal{G}}^2 = \sup_{g \in \mathcal{G}} \text{Var}(g(X))$ and $A_1, A_2$ are universal constants.*

PROOF OF THEOREM 1. We first note that

$$\mathbb{E}\left( K^2\left( \frac{x - X}{h^{1/d}} \right) \right) = h \int_{\mathbb{R}^d} h^{-1} K^2\left( \frac{x - s}{h^{1/d}} \right) f(s)\, ds$$

$$= h \int_{\mathbb{R}^d} K^2(u) f(x - uh^{1/d})\, du \leq h\|f\|_{\infty} \|K\|_2^2,$$

where as usual $\|K\|_2 = (\int_{\mathbb{R}^d} K^2(s)\, ds)^{1/2}$.

Set for $j, k \geq 0$ and $c > 0$, $n_k = 2^k$, $h_{j,k} = (2^j c \log n_k)/n_k$ and

$$\mathcal{K}_{j,k} = \{ K((x - \cdot)/h^{1/d}) : h_{j,k} \leq h \leq h_{j+1,k}, x \in \mathbb{R}^d \}.$$

Clearly for $h_{j,k} \leq h < h_{j+1,k}$, with $\kappa$ as in (K.ii),

$$\mathbb{E}(K^2((x - X)/h^{1/d})) \leq \kappa^2 \wedge 2h_{j,k}\|f\|_{\infty} \|K\|_2^2 =: \kappa^2 \wedge D_0 h_{j,k} =: \sigma_{j,k}^2.$$

We now use Corollary 4 to bound

$$\mathbb{E}\left\|\sum_{i=1}^{n_k}\varepsilon_i g(X_i)\right\|_{\mathcal{K}_{j,k}}.$$

To that end we note that each $\mathcal{K}_{j,k}$ satisfies (i) with $G = \beta = \kappa$. Further, since $\mathcal{K}_{j,k} \subset \mathcal{K}$, we see by (K.iii) that each $\mathcal{K}_{j,k}$ also fulfills (ii). [W.l.o.g. we assume that $\nu, C \geq 1$ in (K.iii).] Noting that

$$C_1\beta/\sigma_0 \leq (\beta^2/\sigma_0^2) \vee C_1^2$$

and the function $h \to h\log(h^{-1} \vee C_1^2)$ is increasing for $h \geq 0$ (recall that $C_1 \geq e$), we see by applying Corollary 4 with $U = \beta = \kappa$ and using the bound $\sigma_0 \leq \sigma_{j,k}$, that we have for $j \geq 0$,

$$\mathbb{E}\left\|\sum_{i=1}^{n_k}\varepsilon_i g(X_i)\right\|_{\mathcal{K}_{j,k}} \leq A\beta\sqrt{\frac{\nu n_k D_0 h_{j,k}}{\beta^2}\log\left(\frac{\beta^2}{D_0 h_{j,k}} \vee C_1^2\right)} + 2A\nu\kappa\log(C_3 n_k),$$

which for $D_1 = A\sqrt{\nu D_0}$ and $D_2 = D_0/\beta^2$ is equal to

$$(2.6) \qquad D_1\sqrt{n_k h_{j,k}\log\left(\frac{1}{D_2 h_{j,k}} \vee C_1^2\right)} + 2A\nu\kappa\log(C_3 n_k).$$

Using once more the fact that $h \to h\log(h^{-1} \vee C_1^2)$ is increasing for $h \geq 0$, we see that the first term of the above bound is, for large $k$, greater than or equal to

$$D_1\sqrt{c\log n_k}\sqrt{\log(n_k/\{cD_2\log n_k\})}.$$

Thus the order of the second term is always smaller than or equal to that of the first one. Consequently, we have for $j \geq 0$ and large enough $k$,

$$\mathbb{E}\left\|\sum_{i=1}^{n_k}\varepsilon_i g(X_i)\right\|_{\mathcal{K}_{j,k}} \leq D_3\sqrt{n_k h_{j,k}\left(\log\left(\frac{1}{D_2 h_{j,k}}\right) \vee \log\log n_k\right)}$$

$$=: D_3 a_{j,k},$$

where $D_3$ is a positive constant.

Applying the Inequality with $M = \kappa$ and $\sigma_{\mathcal{G}}^2 = \sigma_{\mathcal{K}_{j,k}}^2 \leq D_0 h_{j,k}$, we get for any $t > 0$,

$$(2.7) \qquad \begin{aligned} &\mathbb{P}\left\{\max_{n_{k-1} \leq n \leq n_k}\|\sqrt{n}\alpha_n\|_{\mathcal{K}_{j,k}} \geq A_1(D_3 a_{j,k} + t)\right\} \\ &\qquad \leq 2[\exp(-A_2 t^2/(D_0 n_k h_{j,k})) + \exp(-A_2 t/\kappa)]. \end{aligned}$$

Setting for any $\rho > 1$, $j \geq 0$ and $k \geq 1$,

$$p_{j,k}(\rho) = \mathbb{P}\left\{\max_{n_{k-1} \leq n \leq n_k}\|\sqrt{n}\alpha_n\|_{\mathcal{K}_{j,k}} \geq A_1(D_3 + \rho)a_{j,k}\right\},$$

and using the fact that $a_{j,k}^2/n_k h_{j,k} \geq \log \log n_k$, we can infer that for large $k$,

$$p_{j,k}(\rho) \leq 2 \left[ \exp\left( -\frac{\rho^2 A_2}{D_0} \log \log n_k \right) + \exp\left( -\frac{A_2 \rho}{\kappa} \sqrt{n_k h_{j,k} \log \log n_k} \right) \right].$$

Recalling that $h_{j,k} \geq c \log n_k / n_k$, we readily obtain that for large $k$ and any $j \geq 0$,

$$(2.8) \qquad\qquad p_{j,k}(\rho) \leq 4 (\log n_k)^{-\tilde{\rho}},$$

where $\tilde{\rho} = \frac{A_2}{D_0} \rho^2$.

Set $l_k = \max\{j : h_{j,k} \leq 2\}$. It is easy to see that for large $k$,

$$(2.9) \qquad\qquad l_k \leq 2 \log n_k.$$

Hence in view of (2.8) and (2.9) we have for large $k$ and $\rho \geq 1$,

$$P_k(\rho) := \sum_{j=0}^{l_k-1} p_{j,k}(\rho) \leq 8 (\log n_k)^{1-\tilde{\rho}},$$

which implies that if we choose $\rho \geq 2(D_0/A_2)^{1/2}$ (say), we have

$$(2.10) \qquad\qquad \sum_{k=1}^{\infty} P_k(\rho) < \infty.$$

Notice that by definition of $l_k$ for large $k$,

$$2 h_{l_k,k} = h_{l_k+1,k} \geq 2.$$

Consequently, we then have for $n_{k-1} \leq n \leq n_k$,

$$\left[ \frac{c \log n}{n}, 1 \right] \subset \left[ \frac{c \log n_k}{n_k}, h_{l_k,k} \right].$$

Thus for all large enough $k$ and $n_{k-1} \leq n \leq n_k$,

$$A_k(\rho) := \left\{ \max_{n_{k-1} \leq n \leq n_k} \sup_{c \log n/n \leq h \leq 1} \frac{\sqrt{nh} \| \widehat{f}_{n,h} - \mathbb{E} \widehat{f}_{n,h} \|_\infty}{\sqrt{\log(1/h) \vee \log \log n}} > 2 A_1 (D_3 + \rho) \right\}$$

$$\subset \bigcup_{j=0}^{l_k-1} \left\{ \max_{n_{k-1} \leq n \leq n_k} \| \sqrt{n} \alpha_n \|_{\mathcal{K}_{j,k}} \geq A_1 (D_3 + \rho) a_{j,k} \right\},$$

and we see that $\mathbb{P}(A_k(\rho)) \leq P_k(\rho)$. Recalling (2.10), we obtain our theorem via the Borel–Cantelli lemma. $\quad\square$

REMARK 9. We note that if the density is bounded only over $J := I^\varepsilon$, for some $\varepsilon > 0$, with $I$ a compact subset of $\mathbb{R}^d$, and if $K$ is a kernel with support in $[-1/2, 1/2]^d$ we still have for any $0 < h_0 < (2\varepsilon)^d$, with probability 1,

$$(2.11) \qquad \limsup_{n \to \infty} \sup_{c \log n/n \leq h \leq h_0} \frac{\sqrt{nh}\|\widehat{f}_{n,h} - \mathbb{E}\widehat{f}_{n,h}\|_I}{\sqrt{\log(1/h) \vee \log\log n}} =: \tilde{K}(I, c) < \infty.$$

This follows immediately from the above proof by an obvious modification of the bound for $\mathbb{E}K^2((x-X)/h^{1/d})$ and replacing the set $\mathbb{R}^d$ by $I$ in the definition of $\mathcal{K}_{j,k}$.

PROOF OF COROLLARY 1. The proof of (1.3) is obvious from (1.2). Turning to the proof of (1.6), we note that by integrability the assumption that $f$ is uniformly continuous on $\mathbb{R}^d$ is equivalent to $f$ being continuous on $\mathbb{R}^d$ and satisfying the condition that

$$\lim_{R \to \infty} \sup\{f(z) : |z| \geq R\} = 0,$$

which of course implies that $\|f\|_\infty < \infty$. This, when combined with the corollary on page 65 of [28], gives the following lemma.

LEMMA 1. Let $f$ be a uniformly continuous Lebesgue density function on $\mathbb{R}^d$. Then for any kernel $K$ which satisfies (K.i), (K.ii) and (K.v), we have

$$(2.12) \qquad \sup_{z \in \mathbb{R}^d} |f * K_h(z) - f(z)| \to 0 \qquad \text{as } h \searrow 0,$$

where $f * K_h(z) := h^{-1} \int_{\mathbb{R}^d} f(x) K(h^{-1/d}(z-x)) \, dx$.

Observing that $\mathbb{E}\widehat{f}_{n,h}(z) = f * K_h(z)$, we see that Lemma 1 and (1.5) imply (1.6). $\square$

## 3. Proofs of Theorems 2 and 3 and Corollaries 2 and 3.
We are now ready to prove Theorems 2 and 3 and Corollaries 2 and 3. We shall consider a slightly more general setting than in the Introduction, allowing the variables $Y, Y_1, Y_2, \ldots$ to be $r$-dimensional, where $r \geq 1$. Further introduce the following process:

Let $\Phi$ denote a class of measurable functions on $\mathbb{R}^r$ with a finite-valued measurable envelope function $F$, that is,

$$(3.1) \qquad F(y) \geq \sup_{\varphi \in \Phi} |\varphi(y)|, \qquad y \in \mathbb{R}^r.$$

Further assume that $\Phi$ satisfies (K.iii) and (K.iv) with $\mathcal{K}$ replaced by $\Phi$. For any $\varphi \in \Phi$ and continuous functions $c_\varphi$ and $d_\varphi$ on a compact subset of $J$ of $\mathbb{R}^d$, set for $x \in J$,

$$\omega_{\varphi,n,h}(x) = \sum_{i=1}^n (c_\varphi(x)\varphi(Y_i) + d_\varphi(x))K\left(\frac{x - X_i}{h^{1/d}}\right),$$

where $K$ is a kernel with support contained in $[-1/2, 1/2]^d$ such that

$$\sup_{x \in \mathbb{R}^d} |K(x)| =: \kappa < \infty \quad \text{and} \quad \int_{\mathbb{R}^d} K(s)\,ds = 1.$$

For future use introduce two classes of continuous functions on a compact subset $J$ of $\mathbb{R}^d$ indexed by $\Phi$,

$$\mathcal{C} := \{c_\varphi : \varphi \in \Phi\} \quad \text{and} \quad \mathcal{D} := \{d_\varphi : \varphi \in \Phi\}.$$

We shall always assume that the classes $\mathcal{C}$ and $\mathcal{D}$ are relatively compact with respect to the sup-norm topology, which by the Arzela–Ascoli theorem is equivalent to these classes being uniformly bounded and uniformly equicontinuous.

THEOREM 4. *Let $I$ be a compact subset of $\mathbb{R}^d$. Assume that $\Phi$ and $\mathcal{K}$ satisfy the above conditions and the classes of continuous functions $\mathcal{C}$ and $\mathcal{D}$ are as above, that is, relatively compact with respect to the sup-norm topology, where $J = I^\eta$, for some $0 < \eta < 1$. Also assume that*

(3.2) $\qquad\qquad$ *$f$ is continuous and strictly positive on $J$.*

*Further assume that the envelope function $F$ of the class $\Phi$ satisfies*

(3.3) $\qquad\qquad \exists M > 0, \qquad F(Y)\mathbb{1}\{X \in J\} \leq M \qquad a.s.,$

*or for some $p > 2$,*

(3.4) $\qquad\qquad \alpha := \sup_{z \in J} \mathbb{E}(F^p(Y)|X = z) < \infty.$

*Then we have for any $c > 0$ and $0 < h_0 < (2\eta)^d$, with probability 1,*

(3.5) $\quad \limsup_{n \to \infty} \sup_{c(\log n/n)^\gamma \leq h \leq h_0} \frac{\sup_{\varphi \in \Phi} \|\omega_{\varphi,n,h} - \mathbb{E}\omega_{\varphi,n,h}\|_I}{\sqrt{nh(\log(1/h) \vee \log\log n)}} =: P(c) < \infty,$

*where $\gamma = 1$ in the bounded case and $\gamma = 1 - 2/p$ under assumption (3.4).*

Before proving Theorem 4 we shall show how it implies Theorems 2 and 3 and Corollaries 2 and 3. We need the following lemma.

Lemma 2. *Let $\mathcal{H}$ be a class of uniformly equicontinuous functions $g : J \to \mathbb{R}$ and let $K : \mathbb{R}^d \to \mathbb{R}$ be a kernel with support in $[-1/2, 1/2]^d$ so that $\int_{\mathbb{R}^d} K(u)\,du = 1$. Then we have for any sequence of positive constants $b_n \to 0$,*

$$\sup_{g \in \mathcal{H}} \sup_{0 < h < b_n} \|g * K_h - g\|_I \to 0.$$

Proof. A simple transformation shows that if $x \in I$,

$$|g(x) - g * K_h(x)| = \left| \int_{\mathbb{R}^d} (g(x) - g(x - uh^{1/d})) K(u)\,du \right|,$$

which for $h \leq b_n$ and all large enough $n$ is obviously bounded above by

$$\sup\{|g(x) - g(y)| : x, y \in J, |x - y| \leq b_n^{1/d}/2\} \int_{\mathbb{R}^d} |K(u)|\,du.$$

Since the function class $\mathcal{H}$ is uniformly equicontinuous, we readily obtain the assertion of the lemma. $\square$

Proof of Theorem 2. Set

$$\hat{r}_{n,h}(x) = \frac{1}{nh} \sum_{i=1}^{n} Y_i K((x - X_i)/h^{1/d}), \qquad x \in I.$$

Then we obviously have

$$
\begin{aligned}
(3.6) \quad & |\hat{m}_{n,h}(x) - \bar{r}(x,h)/\bar{f}(x,h)| \\
& \leq \frac{1}{|\hat{f}_{n,h}(x)|} |\hat{r}_{n,h}(x) - \bar{r}(x,h)| + \frac{|\bar{r}(x,h)|}{|\hat{f}_{n,h}(x)\bar{f}(x,h)|} |\hat{f}_{n,h}(x) - \bar{f}(x,h)|.
\end{aligned}
$$

From Theorem 4 [setting $r = 1$, $\Phi = \{\varphi_1\}$, where $\varphi_1(y) = y, y \in \mathbb{R}$] it now follows that with probability 1,

$$(3.7) \quad \limsup_{n \to \infty} \sup_{(c \log n/n)^\gamma \leq h \leq b_n} \frac{\sqrt{nh}\|\hat{r}_{n,h}(\cdot) - \bar{r}(\cdot, h)\|_I}{\sqrt{\log(1/h) \vee \log \log n}} < \infty$$

and by (2.11) that

$$(3.8) \quad \limsup_{n \to \infty} \sup_{c \log n/n \leq h \leq b_n} \frac{\sqrt{nh}\|\hat{f}_{n,h}(\cdot) - \bar{f}(\cdot, h)\|_I}{\sqrt{\log(1/h) \vee \log \log n}} < \infty.$$

This last bound of course implies that as $n \to \infty$,

$$\sup_{c \log n/n < h < b_n} \|\hat{f}_{n,h}(\cdot) - \bar{f}(\cdot, h)\|_I = O(1) \qquad \text{a.s.},$$

where we can make the constant in the $O(1)$-term arbitrarily small by choosing $c$ large enough. Combining this observation with the subsequent result following from Lemma 2 that

$$(3.9) \quad \sup_{0 < h \leq b_n} \|\bar{f}(\cdot, h) - f(\cdot)\|_I \to 0$$

and the assumption that the density $f$ is positive on $J$, we can conclude that for $c$ large enough $\hat{f}_{n,h}$ is bounded away from 0 on $I$, uniformly in $c \log n / n < h < b_n$.

Combining this with (1.8) or (1.10) it follows that $\sup_{0 < h \le b_n} \|\bar{r}(\cdot, h)/\bar{f}(\cdot, h)\|_I$ remains bounded. Therefore, we can infer Theorem 2 from (3.6)–(3.8). □

PROOF OF COROLLARY 2. We first note that assumption (1.12) in conjunction with Scheffé's lemma and also condition (1.10) in the unbounded case implies that $r(x) = \mathbb{E}(Y|X = x)f(x)$ is continuous on $J$. Applying Lemma 2 with $\mathcal{H} = \{r\}$, we see that

$$(3.10) \qquad \sup_{0 < h \le b_n} \|\bar{r}(\cdot, h) - r(\cdot)\|_I \to 0,$$

which with (3.9) and Theorem 2 completes the proof of the corollary. □

PROOFS OF THEOREM 3 AND COROLLARY 3. To see how Theorem 3 follows from Theorem 4, set $\Phi = \{\varphi_t\}$, where $\varphi_t(y) = \mathbb{1}\{y \le t\}$, $t, y \in \mathbb{R}$,

$$\mathcal{C} = \{1/f(\cdot)\} \quad \text{and} \quad \mathcal{D} = \{-F(t|\cdot)/f(\cdot) : t \in \mathbb{R}\}.$$

The classes $\Phi$ and $\mathcal{C}$ clearly satisfy the assumptions of Theorem 4. To see that the function class $\mathcal{D}$ is a relatively compact class of continuous functions on $J$ refer to pages 6 and 7 of [11], which also implies that the class

$$\mathcal{H} = \{g_t : t \in \mathbb{R}\},$$

where for each $t \in \mathbb{R}$, $g_t(\cdot) = F(t|\cdot)f(\cdot)$, is also a relatively compact class of functions defined on $J$. Therefore Theorem 3 and Corollary 3 follow in the same way that Theorem 2 and Corollary 2 did from Theorem 4 and Lemma 2. □

PROOF OF THEOREM 4. We first note that

$$\limsup_{n \to \infty} \sup_{c \log n / n \le h \le h_0} \frac{\sup_{\varphi \in \Phi} \|d_\varphi(x) \sum_{i=1}^n \{K(\frac{x - X_i}{h^{1/d}}) - \mathbb{E}K(\frac{x - X_i}{h^{1/d}})\}\|_I}{\sqrt{nh(\log(1/h) \vee \log \log n)}}$$

$$\le \sup_{\varphi \in \Phi} \|d_\varphi\|_I \limsup_{n \to \infty} \sup_{c \log n / n \le h \le h_0} \frac{\|\sum_{i=1}^n \{K(\frac{x - X_i}{h^{1/d}}) - \mathbb{E}K(\frac{x - X_i}{h^{1/d}})\}\|_I}{\sqrt{nh(\log(1/h) \vee \log \log n)}}.$$

In view of (2.11) it is obvious that this quantity is finite with probability 1. Therefore, if we set for $\varphi \in \Phi$, $x \in I$ and $h > 0$,

$$(3.11) \qquad \eta_{\varphi, n, h}(x) = c_\varphi(x) \sum_{i=1}^n \varphi(Y_i) K\left(\frac{x - X_i}{h^{1/d}}\right),$$

it clearly suffices to show:

PROPOSITION 2. *Under the assumptions of Theorem* 4, *for all* $c > 0$, *there exists a* $Q(c) > 0$ *such that with probability* 1,

$$(3.12) \quad \limsup_{n \to \infty} \sup_{c(\log n/n)^\gamma \leq h \leq h_0} \frac{\sup_{\varphi \in \Phi} \|\eta_{\varphi,n,h} - \mathbb{E}\eta_{\varphi,n,h}\|_I}{\sqrt{nh(\log(1/h) \vee \log\log n)}} =: Q(c).$$

PROOF. We shall prove Proposition 2 under assumption (3.4), as it follows in the bounded case directly from the proof of Theorem 1 and Remark 9. Just replace the classes $\mathcal{K}_{j,k}$ by the classes

$$\mathcal{G}_{j,k} = \{(y,z) \to \varphi(y)c_\varphi(z)K((x - z)/h^{1/d}) : \varphi \in \Phi, x \in I, h_{j,k} \leq h \leq h_{j+1,k}\}.$$

Observe that these classes also satisfy conditions (i)–(iii) of Proposition 1 with $G = \beta = \kappa \sup_{\varphi \in \Phi} \|c_\varphi\|_I$ and our proof of Theorem 1 still works after some minor modifications.

We turn to the unbounded case, that is, assume (3.4) for some $p > 2$. Recall that $\gamma = 1 - 2/p$. For any $k = 1, 2, \ldots$ and $\varphi \in \Phi$, set $n_k = 2^k, a_k = c(\log n_k/n_k)^\gamma$ and

$$(3.13) \quad \varphi_k(y) = \varphi(y)\mathbb{1}\{F(y) < (n_k/k)^{1/p}\}.$$

For $n_{k-1} \leq n \leq n_k$, $x \in I$, $a_k \leq h \leq h_0$ and $\varphi \in \Phi$, let

$$(3.14) \quad \eta_{\varphi,n,h}^{(k)}(x) = c_\varphi(x) \sum_{i=1}^n \varphi_k(Y_i)K\left(\frac{x - X_i}{h^{1/d}}\right).$$

The proof of Proposition 2 in the unbounded case will be a consequence of two lemmas. We will first show:

LEMMA 3. *There exists a constant* $Q_1(c) < \infty$, *such that with probability* 1,

$$(3.15) \quad \limsup_{k \to \infty} \Delta_k = Q_1(c),$$

*where*

$$\Delta_k = \max_{n_{k-1} \leq n \leq n_k} \sup_{a_k \leq h \leq h_0} \frac{\sup_{\varphi \in \Phi} \|\eta_{\varphi,n,h}^{(k)} - \mathbb{E}\eta_{\varphi,n,h}^{(k)}\|_I}{\sqrt{nh(\log(1/h) \vee \log\log n)}}.$$

PROOF. For $x \in I$, $a_k \leq h \leq h_0$ and $\varphi \in \Phi$, let

$$(3.16) \quad v_{\varphi,h,x}(u,v) = c_\varphi(x)\varphi(v)K\left(\frac{x - u}{h}\right)$$

and

$$(3.17) \quad v_{\varphi,h,x}^{(k)}(u,v) = c_\varphi(x)\varphi_k(v)K\left(\frac{x - u}{h}\right).$$

Notice that

$$\eta_{\varphi,n,h}^{(k)}(x) - \mathbb{E}\eta_{\varphi,n,h}^{(k)} = n^{1/2}\alpha_n(v_{\varphi,h,x}^{(k)}), \tag{3.18}$$

where $\alpha_n$ is the empirical process based on $(X_1, Y_1), \ldots, (X_n, Y_n)$. For $k \geq 1$ let

$$\mathcal{G}_k(h) := \{v_{\varphi,h,x}^{(k)} : \varphi \in \Phi \text{ and } x \in I\}.$$

We see that

$$\Delta_k = \max_{n_{k-1} \leq n \leq n_k} \sup_{a_k \leq h \leq h_0} \frac{\|\alpha_n\|_{\mathcal{G}_k(h)}}{\sqrt{nh(\log(1/h) \vee \log\log n)}}.$$

Note that for each $v_{\varphi,h,x}^{(k)} \in \mathcal{G}_k(h)$,

$$\|v_{\varphi,h,x}^{(k)}\|_\infty \leq \|K\|_\infty \sup_{\varphi \in \Phi}\|c_\varphi\|_\infty (n_k/k)^{1/p} =: D_0(n_k/k)^{1/p}. \tag{3.19}$$

Also observe that

$$\mathbb{E}[(v_{\varphi,h,x}^{(k)})^2(X,Y)] \leq \mathbb{E}[v_{\varphi,h,x}^2(X,Y)] \leq \mathbb{E}\left[(c_\varphi(x)\varphi(Y))^2 K^2\left(\frac{x-X}{h^{1/d}}\right)\right].$$

Using a conditioning argument, we infer that this last term is uniformly over $x \in I$

$$\leq \|c_\varphi\|_J^2 \int_{|x-t| \leq h^{1/d}/2} \mathbb{E}[F^2(Y)|X=t]f_X(t)K^2\left(\frac{x-t}{h^{1/d}}\right)dt$$

$$\leq \|c_\varphi\|_J^2 \alpha^{2/p} \int_{|x-t| \leq h^{1/d}/2} f_X(t)K^2\left(\frac{x-t}{h^{1/d}}\right)dt$$

$$\leq h\|c_\varphi\|_J^2 \alpha^{2/p} \int_{[-1/2,1/2]^d} f_X(x-h^{1/d}u)K^2(u)\,du$$

$$\leq h\alpha^{2/p} \sup_{\varphi \in \Phi}\|c_\varphi\|_J^2 \|f_X\|_J\|K\|_2^2 =: hD_1.$$

Thus

$$\sup_{v \in \mathcal{G}_k(h)} \mathbb{E}v^2(X,Y) \leq hD_1. \tag{3.20}$$

Set for $j, k \geq 0$, $h_{j,k} = 2^j a_k$ and

$$\mathcal{G}_{j,k} := \{v_{\varphi,h,x}^{(k)} : \varphi \in \Phi, \ x \in I \text{ and } h_{j,k} \leq h \leq h_{j+1,k}\}.$$

Clearly by (3.20), for all $x \in I$, $\varphi \in \Phi$ and $h_{j,k} \leq h < h_{j+1,k}$,

$$\mathbb{E}[(v_{\varphi,h,x}^{(k)})^2(X,Y)] \leq hD_1 \leq 2D_1h_{j,k},$$

which gives

$$(3.21) \qquad \sup_{v \in \mathcal{G}_{j,k}} \mathbb{E}v^2(X,Y) \leq 2D_1 h_{j,k}.$$

We shall use Corollary 4 to bound $\mathbb{E}\left\|\sum_{i=1}^{n_k} \varepsilon_i v(X_i, Y_i)\right\|_{\mathcal{G}_{j,k}}$. Note first that by arguing as in the proof of Lemma 5 of EM [11], each $\mathcal{G}_{j,k} \subset \mathcal{G}$, where $\mathcal{G}$ is a class that satisfies (K.iii) and (K.iv) with $\mathcal{K}$ replaced by $\mathcal{G}$. Next by setting $U = D_0(n_k/k)^{1/p}$ and $\beta^2 = \alpha^{2/p}$ it follows that

$$\mathbb{E}\left\|\sum_{i=1}^{n_k} \varepsilon_i v(X_i, Y_i)\right\|_{\mathcal{G}_{j,k}} \leq A\sqrt{\nu n_k \sigma_0^2 \log(C_1 \beta/\sigma_0)}$$
$$+ 2AD_0 \nu (n_k/k)^{1/p} \log(C_3 D_0^{-2} \beta^2 n_k^{1-2/p} k^{2/p}).$$

Replacing $\sigma_0^2$ in the first term by the upper bound $2D_1 h_{j,k} \wedge \beta^2$ and recalling that

$$h_{j,k} \geq a_k = c(\log n_k/n_k)^{1-2/p},$$

we obtain after a small calculation that for suitable positive constants $D_2$ and $D_3$,

$$(3.22) \qquad \mathbb{E}\left\|\sum_{i=1}^{n_k} \varepsilon_i v(X_i, Y_i)\right\|_{\mathcal{G}_{j,k}} \leq D_3 \sqrt{n_k h_{j,k} \log((D_2 h_{j,k})^{-1} \vee C_1)}.$$

Set

$$a_{j,k} = \sqrt{n_k h_{j,k}\left(\log\left(\frac{1}{D_2 h_{j,k}}\right) \vee \log\log n_k\right)}, \qquad k \geq 1, j \geq 0.$$

Applying Talagrand's inequality with $M = D_0(n_k/k)^{1/p}$ and $\sigma_{\mathcal{G}}^2 = \sigma_{\mathcal{G}_{j,k}}^2 \leq 2D_1 h_{j,k}$, we get for any $t > 0$ and large enough $k$,

$$\mathbb{P}\left\{\max_{n_{k-1} \leq n \leq n_k} \|\sqrt{n}\alpha_n\|_{\mathcal{G}_{j,k}} \geq A_1(D_3 a_{j,k} + t)\right\}$$
$$\leq 2[\exp(-A_2 t^2/(2D_1 n_k h_{j,k})) + \exp(-A_2 t k^{1/p}/(D_0 n_k^{1/p}))].$$

Set for any $\rho > 1$, $j \geq 0$ and $k \geq 1$,

$$p_{j,k}(\rho) = \mathbb{P}\left\{\max_{n_{k-1} \leq n \leq n_k} \|\sqrt{n}\alpha_n\|_{\mathcal{G}_{j,k}} \geq A_1(D_3 + \rho)a_{j,k}\right\}.$$

Using the facts that $a_{j,k}^2/(n_k h_{j,k}) \geq \log\log n_k$ and $h_{j,k} \geq c(\log n_k/n_k)^{1-2/p}$, we readily obtain for large $k$ and $j \geq 0$ that

$$p_{j,k}(\rho) \leq 2\exp\left(-\frac{\rho^2 A_2}{D_1} \log\log n_k\right) + 2\exp\left(-\frac{\sqrt{c}\rho A_2}{D_0}\sqrt{\log n_k \log\log n_k}\right),$$

which for $\tilde{\rho} = \frac{A_2}{D_1}\rho^2$ and large $k$ is less than or equal to $4(\log n_k)^{-\tilde{\rho}}$. Set $l_k = \max\{j : h_{j,k} \leq 2h_0\}$ if this set is nonempty, which is obviously the case for large enough $k$. Then we have $l_k \leq k$ for large enough $k$ and, consequently,

$$P_k(\rho) := \sum_{j=0}^{\ell_k - 1} p_{j,k}(\rho) \leq 4\ell_k(\log n_k)^{-\tilde{\rho}} \leq k^{-2},$$

provided we have chosen $\rho$ large enough.

Further notice that by the definition of $l_k$ for large $k$,

$$2h_{\ell_k,k} = h_{l_k+1,k} \geq 2h_0,$$

which implies that we have for $n_{k-1} \leq n \leq n_k$,

$$[a_k, h_0] \subset [a_k, h_{\ell_k,k}].$$

Thus for all large enough $k$,

$$A_k(\rho) := \{\Delta_k \geq 2A_1(D_3 + \rho)\} \subset \bigcup_{j=0}^{\ell_k-1}\left\{\max_{n_{k-1} \leq n \leq n_k} \|\sqrt{n}\alpha_n\|_{\mathcal{G}_{j,k}} \geq A_1(D_3 + \rho)a_{j,k}\right\}.$$

It follows now with the above choice for $\rho$,

$$\mathbb{P}(A_k(\rho)) \leq P_k(\rho) \leq k^{-2},$$

which by the Borel–Cantelli lemma implies Lemma 3.  $\square$

Write

$$(3.23) \qquad \overline{\varphi}_k(y) = \varphi(y)\mathbb{1}\{F(y) \geq (n_k/k)^{1/p}\}.$$

For $\varphi \in \Phi$, $x \in I$ and $n_{k-1} \leq n \leq n_k$, let

$$(3.24) \qquad \overline{\eta}^{(k)}_{\varphi,n,h}(x) = c_\varphi(x)\sum_{i=1}^{n}\left\{\overline{\varphi}_k(Y_i)K\left(\frac{x - X_i}{h^{1/d}}\right)\right\}.$$

LEMMA 4.  *With probability* 1,

$$(3.25) \quad \lim_{k \to \infty}\max_{n_{k-1} \leq n \leq n_k}\sup_{c(\log n_k/n_k)^{1-2/p} \leq h \leq h_0}\frac{\sup_{\varphi \in \Phi}\|\overline{\eta}^{(k)}_{\varphi,n,h} - \mathbb{E}\overline{\eta}^{(k)}_{\varphi,n,h}\|_I}{\sqrt{nh(\log(1/h) \vee \log\log n)}} = 0.$$

PROOF.   First note that for any $h \leq h_0, \varphi \in \Phi$ and $n_{k-1} \leq n \leq n_k$,

$$\|\mathbb{E}\overline{\eta}^{(k)}_{\varphi,n,h}\|_I \leq \kappa \sup_{\varphi \in \Phi}\|c_\varphi\|_I n_k \mathbb{E}[F(Y)\mathbb{1}\{X \in J, F(Y) \geq (n_k/k)^{1/p}\}].$$

We further have by (3.4),

$$\mathbb{E}F^p(Y)\mathbb{1}\{X \in J\} < \infty,$$

and we see that uniformly in $n_{k-1} \leq n \leq n_k$, $h \leq h_0$ and $\varphi \in \Phi$,

$$\|\mathbb{E}\overline{\eta}_{\varphi,n,h}^{(k)}\|_I = o(n_k^{1/p} k^{1-1/p}) = o(\sqrt{n_k a_k \log(1/a_k)})$$

as $k \to \infty$, where $a_k = c(\log n_k / n_k)^{1-2/p}$.

By monotonicity of the function $h \to h \log(1/h), h \leq 1/e$, we readily obtain that

$$(3.26) \quad \lim_{k \to \infty} \max_{n_{k-1} \leq n \leq n_k} \sup_{a_k \leq h \leq h_0} \frac{\sup_{\varphi \in \Phi} \|\mathbb{E}\overline{\eta}_{\varphi,n,h}^{(k)}\|_I}{\sqrt{nh(\log(1/h) \vee \log\log n)}} = 0.$$

It remains to be shown that, with probability 1,

$$(3.27) \quad \lim_{k \to \infty} \max_{n_{k-1} \leq n \leq n_k} \sup_{a_k \leq h \leq h_0} \frac{\sup_{\varphi \in \Phi} \|\overline{\eta}_{\varphi,n,h}^{(k)}\|_I}{\sqrt{nh(\log(1/h) \vee \log\log n)}} = 0.$$

Similarly as above we have

$$\max_{n_{k-1} \leq n \leq n_k} \sup_{a_k \leq h \leq h_0} \sup_{\varphi \in \Phi} \|\overline{\eta}_{\varphi,n,h}^{(k)}\|_I$$

$$\leq \kappa \sup_{\varphi \in \Phi} \|c_\varphi\|_I \sum_{i=1}^{n_k} F(Y_i) \mathbb{1}\{X_i \in J, F(Y_i) > (n_k/k)^{1/p}\}.$$

Inspecting the proof of Lemma 1 of [11], we see that the argument there also applies if we set $h_n = c(\log n/n)^{1-2/p}$ to give

$$\sum_{i=1}^{n_k} F(Y_i) \mathbb{1}\{X_i \in J, F(Y_i) > (n_k/k)^{1/p}\} = o(\sqrt{n_k a_k \log(1/a_k)}),$$

as $k \to \infty$, and we see by the same argument as in (3.26) that (3.27) holds, thereby finishing the proof of Lemma 4. $\square$

Proposition 2 now follows from Lemmas 3 and 4. $\square$

**Acknowledgments.** The authors thank Luc Devroye and Jon Wellner for posing the questions that led to our work.

DEPARTEMENT WISKUNDE
VRIJE UNIVERSITEIT BRUSSEL
PLEINLAAN 2
B-1050 BRUSSEL
BELGIUM
E-MAIL: ueinmahl@vub.ac.be

DEPARTMENT OF FOOD
  AND RESOURCE ECONOMICS
UNIVERSITY OF DELAWARE
206 TOWNSEND HALL
NEWARK, DELAWARE 19717
USA
E-MAIL: davidm@udel.edu